%
%
\input amstex.tex
\documentstyle{amsppt}
\magnification=1200
\baselineskip=13pt
\hsize=6.5truein
\vsize=8.9truein
\countdef\sectionno=1
\countdef\eqnumber=10
\countdef\theoremno=11
\countdef\countrefno=12
\countdef\cntsubsecno=13
\sectionno=0
\def\newsection{\global\advance\sectionno by 1
                \global\eqnumber=1
                \global\theoremno=1
                \global\cntsubsecno=0
                \the\sectionno}

\def\newsubsection#1{\global\advance\cntsubsecno by 1
                     \xdef#1{{\S\the\sectionno.\the\cntsubsecno}}
                     \ \the\sectionno.\the\cntsubsecno.}

\def\theoremname#1{\the\sectionno.\the\theoremno
                   \xdef#1{{\the\sectionno.\the\theoremno}}
                   \global\advance\theoremno by 1}

\def\eqname#1{\the\sectionno.\the\eqnumber
              \xdef#1{{\the\sectionno.\the\eqnumber}}
              \global\advance\eqnumber by 1}

\def\thmref#1{#1}

\global\countrefno=1

\def\refno#1{\xdef#1{{\the\countrefno}}\global\advance\countrefno by 1}

\def\R{{\Bbb R}}
\def\N{{\Bbb N}}
\def\C{{\Bbb C}}
\def\Z{{\Bbb Z}}
\def\Zp{{\Bbb Z}_+}
\def\hf{{1\over 2}}
\def\hZp{{1\over 2}\Zp}
\def\A{{\Cal A}_q(SU(2))}
\def\U{{\Cal U}_q({\frak{su}}(2))}
\def\Hi{\ell^2(\Zp)}
\def\a{\alpha}
\def\b{\beta}
\def\g{\gamma}
\def\d{\delta}
\def\s{\sigma}
\def\t{\tau}
\def\l{\lambda}
\def\m{\mu}
\def\n{\nu}
\def\r{\rho_{\t ,\s}}
\def\ris{\rho_{\infty ,\s}}
\def\rti{\rho_{\t ,\infty}}
\def\p{\pi}
\def\vp{\varphi}
\refno{\AndrAbqJ}
\refno{\Aske}
\refno{\AskeW}
\refno{\Bere}
\refno{\Domb}
\refno{\GaspR}
\refno{\Jimb}
\refno{\KalnM}
\refno{\KoelSIAM}
\refno{\KoelAF}
\refno{\KoelAAM}
\refno{\KoelV}
\refno{\KoorIM}
\refno{\KoorOPTA}
\refno{\KoorAF}
\refno{\KoorZSE}
\refno{\MasuJFA}
\refno{\Noum}
\refno{\NoumMJASJ}
\refno{\NoumMLNM}
\refno{\NoumMCM}
\refno{\Rahm }
\refno{\RahmV}
\refno{\VaksS}
\refno{\VAsscK}
\refno{\WoroRIMS}
\refno{\WoroCMP}

\topmatter
\title Addition formula for $2$-parameter family of Askey-Wilson polynomials
\endtitle
\rightheadtext{Addition formula for Askey-Wilson polynomials}
\author H.T. Koelink\endauthor
\address Departement Wiskunde, Katholieke Universiteit Leuven,
Celestijnenlaan 200 B, B-3001 Leuven (Heverlee), Belgium\endaddress
\email erik.koelink\@wis.kuleuven.ac.be\endemail
\date December 1994\enddate
\thanks Supported by a Fellowship of the Research Council of the
Katholieke Universiteit Leuven. \endthanks
\abstract
For a two parameter family of Askey-Wilson polynomials, that can be
regarded as basic analogues of the Legendre polynomials, an addition
formula is derived. The addition formula is a two-parameter extension of
Koornwinder's addition formula for the little $q$-Legendre polynomials.
A corresponding product formula is derived. As the base tends to one, the
addition and product formula go over into the addition and product
formula for the Legendre polynomial. The addition formula is derived
from the interpretation of Askey-Wilson polynomials as generalised
matrix elements on the quantum $SU(2)$ group.
\endabstract
\endtopmatter
\document

\head\newsection . Introduction \endhead

Let us start with describing the classical addition formula for the
Legendre polynomials.
The Jacobi polynomials $R_n^{(\a,\b)}(x)$ are orthogonal polynomials
with respect to the beta measure $(1-x)^\a (1+x)^\b$ on the interval
$[-1,1]$ and normalised by $R_n^{(\a,\b)}(1)=1$. The polynomial
$R_n^{(0,0)}(x)$ is the Legendre polynomial, and for $\a=\b$ the
name Gegenbauer or ultraspherical polynomials is used. Laplace proved
the following addition formula for the Legendre polynomials in 1782;
$$
\gathered
R_l^{(0,0)}\bigl( xy+t\sqrt{(1-x^2)(1-y^2)}\bigr)
= R^{(0,0)}_l(x)R^{(0,0)}_l(y) \\
+2\sum_{n=1}^l {{(l+n)!}\over{(l-n)!(n!)^2}}2^{-2n}
(1-x^2)^{n/2} R^{(n,n)}_{l-n}(x)
(1-y^2)^{n/2} R^{(n,n)}_{l-n}(y)T_n(t),
\endgathered
\tag\eqname{\vgladdformLeg}
$$
where $T_n(\cos\theta)=\cos n\theta = R_n^{(-1/2,-1/2)}(\cos\theta)$, is the
Chebyshev polynomial of the first kind, cf. Askey \cite{\Aske, Lecture~4}.
The orthogonality relations for the Chebyshev polynomials yield the
product formula
$$
\gathered
R^{(n,n)}_{l-n}(x)R^{(n,n)}_{l-n}(y)=  2^{2n}
{{(l-n)!(n!)^2}\over{\pi (l+n)!}} \bigl( \sqrt{(1-x^2)(1-y^2)}\bigr)^{-n} \\
\times \int_{-1}^1 R^{(0,0)}_l\bigl( xy+t\sqrt{(1-x^2)(1-y^2)}\bigr)
{{T_n(t)}\over{\sqrt{1-t^2}}}\, dt,
\endgathered
\tag\eqname{\vglprodformultra}
$$
which for $n=0$ is the product formula for the Legendre polynomial.

For basic analogues of the Legendre polynomials there are several
analogues of the addition formula \thetag{\vgladdformLeg}. The addition
formula for the continuous $q$-Legendre polynomials is proved by Rahman
and Verma \cite{\RahmV} as a special case of their addition formula for
the continuous $q$-ultraspherical polynomials, which have been
introduced by Rogers in 1895. A quantum $SU(2)$ group theoretic proof of the
addition formula for the continuous $q$-Legendre polynomials has been
given by the author \cite{\KoelSIAM}. Also, the quantum $SU(2)$ group has
naturally led to an addition formula for the little $q$-Legendre
polynomials, cf. Koornwinder \cite{\KoorAF},
since the matrix elements of the irreducible unitary
representations of the quantum $SU(2)$ group are explicitly known in
terms of little $q$-Jacobi polynomials \cite{\VaksS}, \cite{\MasuJFA},
\cite{\KoorIM}. See Rahman \cite{\Rahm} for
an analytic proof of this addition formula. There is also an addition
formula for the big $q$-Legendre polynomials \cite{\KoelAF}, and the
addition formula derived in this paper contains these last two addition
formulas as limiting cases.
Apart from the quantum group theoretic approach it is also possible to use
quantum algebras and its representation theory
in proving addition formulas for $q$-special functions, see
e.g. Kalnins and Miller \cite{\KalnM} and references therein.

As a result of Koornwinder's work \cite{\KoorZSE}
on zonal spherical elements on the quantum $SU(2)$ group and
Askey-Wilson polynomials, there have been
papers by Noumi and Mimachi \cite{\NoumMJASJ}, \cite{\NoumMLNM} and the
author \cite{\KoelSIAM}, \cite{\KoelAAM}, in which the full four
parameter family of Askey-Wilson polynomials is interpreted as
generalised matrix elements on the quantum $SU(2)$ group. This has led
to an abstract addition formula for Askey-Wilson polynomials, i.e.
involving non-commuting variables, \cite{\NoumMJASJ}, \cite{\NoumMLNM},
\cite{\KoelSIAM}, \cite{\KoelAAM}, from which a (degenerate) addition
formula for Askey-Wilson polynomials can be obtained.

In this paper we derive an addition formula for a two-parameter family
of Askey-Wilson polynomials much along the lines of the quantum group
theoretic proof of the addition formula for the little $q$-Legendre
polynomials. The proof is based on the fact that we have a basis of
eigenvectors for a certain self-adjoint operator in an irreducible
$\ast$-representation of the $C^\ast$-algebra for the quantum $SU(2)$
group \cite{\KoelAF}
and that the non-polynomial minimal part of the associated spherical
generalised elements can be factorised in elements which act nicely in
this basis of eigenvectors. This last part is suggested by the paper by
Noumi and Mimachi \cite{\NoumMCM}, and we make an explicit
identification with their work in remark 3.5. 
All this is explained in \S 3. The proof of the addition formula is then
relatively easy and is given in \S 4. In \S 5 the corresponding product
formula is derived and the limit transition $q\uparrow 1$ to
\thetag{\vgladdformLeg} and \thetag{\vglprodformultra}
is considered in \S 6. In \S 2 we recall the necessary results on the
relation between certain basic hypergeometric orthogonal polynomials and
the quantum $SU(2)$ group.

To end  this introduction we remark that this paper is concerned with
transforming an identity for $q$-special functions involving non-commuting
variables into an identity for $q$-special functions in commuting variables.

\head\newsection . Orthogonal polynomials and the quantum $SU(2)$ group
\endhead

In this section we present some of the relations between basic hypergeometric
orthogonal polynomials and the quantum $SU(2)$ group that are necessary to
obtain an addition formula for a two-parameter family of Askey-Wilson
polynomials. More information can be found in the survey papers by Koornwinder
\cite{\KoorOPTA}, Noumi \cite{\Noum} and the author \cite{\KoelAAM}.

$\A$ is the complex unital associative algebra generated
by $\a$, $\b$, $\g$, $\d$ subject to the relations
$$
\gathered
\a\b =q\b\a ,\quad \a\g = q\g\a ,\quad \b\d = q\d\b ,\quad \g\d = q\d\g
,\\
\b\g =\g\b ,\quad \a\d -q\b\g = \d\a - q^{-1}\b\g =1
\endgathered
\tag\eqname{\vglcommrelAq}
$$
for some constant $q\in\C$.
With a $\ast$-operator given by
$$
\a^\ast = \d , \quad \b^\ast = -q\g ,\quad \g^\ast = -q^{-1}\b ,\quad
\d^\ast = \a
\tag\eqname{\vglsteropAq}
$$
the algebra $\A$ becomes a $\ast$-algebra for non-zero real $q$, and
from now on we fix $0<q<1$. The $\ast$-algebra can be completed into a
$C^\ast$-algebra, and then this quantum group has been introduced by
Woronowicz \cite{\WoroRIMS} as an example of his general theory of compact
matrix quantum groups \cite{\WoroCMP}.

All $\ast$-representations of the $C^\ast$-algebra have been
classified, cf. Vaksman and Soibelman \cite{\VaksS}, and we only use the
infinite dimensional $\ast$-representation $\p$ of $\A$ in the Hilbert space
$\Hi$ with orthonormal basis $\{ e_n\mid n\in\Zp\}$ given by
$$
\gathered
\p(\a) e_n = \sqrt{1-q^{2n}} e_{n-1}, \quad
\p(\b )e_n = -q^{n+1} e_n, \\
\p(\g )e_n = q^n e_n, \quad
\p(\d) e_n = \sqrt{1-q^{2n+2}} e_{n+1},
\endgathered
\tag\eqname{\vgldefreprA}
$$
where we use the convention $e_{-p}=0$ for $p\in\N$. Note that
$-q\p(\g)=\p(\b)$. For all $\xi\in\A$, $\p(\xi)$ is a bounded operator
on $\Hi$.

The algebra $\A$ is an example of a Hopf $\ast$-algebra.
The comultiplication $\Delta$, which is a $\ast$-homomorphism of
$\A\to\A\otimes\A$, is given on the generators by
$$
\gathered
\Delta(\a)=\a\otimes\a + \b\otimes\g ,\quad
\Delta(\b )=\a\otimes\b + \b\otimes\d ,\\
\Delta(\g )=\g\otimes\a + \d\otimes\g ,\quad
\Delta(\d )=\g\otimes\b + \d\otimes\d ,
\endgathered
\tag\eqname{\vgldefDeltaonAq}
$$
The representations of the quantum $SU(2)$ group
correspond to the corepresentations of the Hopf $\ast$-algebra $\A$.
The irreducible unitary representations of the
quantum $SU(2)$ group are completely classified and, as for the Lie group
$SU(2)$, they are indexed by the spin $l$, $l\in\hZp$, for each dimension
$2l+1$, cf. Woronowicz
\cite{\WoroRIMS}, \cite{\WoroCMP}, Koornwinder \cite{\KoorIM}, Masuda et al.
\cite{\MasuJFA}, Vaksman and Soibelman \cite{\VaksS}. The matrix elements
$t^l_{n,m}\in\A$, $n,m=-l,-l+1,\ldots,l$, of the irreducible unitary
corepresentation $t^l$
of spin $l$ are explicitly known in terms of the little $q$-Jacobi polynomials
\cite{\KoorIM}, \cite{\MasuJFA}, \cite{\VaksS}. The cohomorphism property for
the matrix elements is
$$
\Delta (t^l_{n,m}) = \sum_{k=-l}^l t^l_{n,k}\otimes t^l_{k,m},
\tag\eqname{\vglcohomtl}
$$
and for $n=m=0$ Koornwinder \cite{\KoorAF} has transformed this identity into
an addition formula for the little $q$-Legendre polynomials.

It is possible to extend the notion of matrix elements to get generalised
matrix elements, which satisfy relative left and right infinitesimal
invariance. In the spherical case this has been introduced by Koornwinder
\cite{\KoorOPTA}, \cite{\KoorZSE} and for the associated spherical elements it
has been worked out by the author \cite{\KoelSIAM} and for generalised matrix
elements Noumi and Mimachi \cite{\NoumMJASJ}, \cite{\NoumMLNM} have stated the
corresponding theorems, see the survey \cite{\KoelAAM} for a complete proof.
To describe the infinitesimal invariance the duality of the Hopf
$\ast$-algebra with the quantised universal enveloping algebra $\U$ as
introduced by Jimbo \cite{\Jimb} is needed, but we only need some of the
results, see \cite{\KoelAAM} for proofs.
There exist elements $b^l_{i,j}(\t,\s)\in\A$, $l\in\hZp$, $i,j\in\{
-l,-l+1,\ldots,l-1,l\}$, $\s,\t\in\R\cup\{ \infty\}$ such that
$$
\Delta\bigl(b^l_{i,j}(\t,\s)\bigr) = \sum_{n=-l}^l
\bigl( D.b^l_{i,n}(\t,\m)\bigr) \otimes b^l_{n,j}(\m,\s),
\qquad\forall\;\m\in\R\cup\{\infty\}.
\tag\eqname{\vglabstractaddf}
$$
Here $D.\colon\A\to\A$ is the bijective algebra homomorphism given
by $D.\a=q^{-\hf}\a$, $D.\b=q^\hf\b$, $D.\g=q^{-\hf}\g$ and $D.\d=q^\hf\d$.
For $\s=\t=\infty$ the
generalised matrix element $b^l_{i,j}(\t,\s)$ reduces to $t^l_{i,j}$ up to a
constant, and $\s=\t=\m=\infty$ in \thetag{\vglabstractaddf} reduces to
\thetag{\vglcohomtl}.
The case of most interest to us is $l\in\Zp$, $i=j=0$ and $\m=\infty$ in
\thetag{\vglabstractaddf}; i.e.
$$
\Delta\bigl(b^l_{0,0}(\t,\s)\bigr) = \sum_{n=-l}^l
\bigl( D.b^l_{0,n}(\t,\infty)\bigr) \otimes b^l_{n,0}(\infty,\s).
\tag\eqname{\vglabstractadfeen}
$$

The generalised matrix elements are explicitly known in terms of Askey-Wilson
polynomials, where the argument is a simple element of the Hopf $\ast$-algebra
$\A$, cf. \cite{\KoorOPTA}, \cite{\KoorZSE}, \cite{\NoumMJASJ},
\cite{\NoumMLNM}, \cite{\KoelSIAM}, \cite{\KoelAAM}. Define in $\A$
$$
\aligned
\r = {1\over 2}\bigl(& \a^2+\d^2 +q\g^2 +q^{-1}\b^2 + i(q^{-\s}-q^\s )
(q\d\g +\b\a )\\
& -i(q^{-\t}-q^\t)(\d\b +q\g\a) + (q^{-\s}-q^\s)(q^{-\t}-q^\t)\b\g
\bigr),
\endaligned
\tag\eqname{\vgldefrhosigmatau}
$$
then
$$
b^l_{0,0}(\t,\s)= q^{-l} (q^{2l+2};q^2)_l^{-1}
p^{(0,0)}_l(\r;q^\t, q^\s \mid q^2).
\tag\eqname{\vglzonalseAW}
$$
The normalisation has been chosen differently from \cite{\KoelAAM}. It is
chosen such that, with \thetag{2.12}, 
the common constant in \thetag{\vglabstractadfeen} is dropped.
Here we use the notation
$$
p_n^{(\a,\b)}(x;s,t\mid q) = p_n(x;q^{1/2+\a}s/t,q^{1/2}t/s,
-q^{1/2}/(st), -stq^{1/2+\b}\mid q)
\tag\eqname{\vgldefAWalsqJacobi}
$$
for $q$-analogues of the Jacobi polynomial, where
$$
p_n(\cos\theta;a,b,c,d\mid q) = a^{-n} (ab,ac,ad;q)_n
\, {}_4\vp_3 \left( {{q^{-n}, abcdq^{n-1}, ae^{i\theta},a^{-i\theta}}
\atop{ab,\ ac,\ ad}};q,q\right)
\tag\eqname{\vgldefAWpol}
$$
denotes an Askey-Wilson polynomial \cite{\AskeW}. The notation for basic
hypergeometric series follows the excellent book \cite{\GaspR} by Gasper and
Rahman;
$$
{}_{r+1}\vp_r \left( {{a_1,\ldots,a_{r+1}}\atop{b_1,\ldots, b_r}};q,z\right)
= \sum_{k=0}^\infty {{(a_1,\ldots,a_{r+1};q)_k}\over{(q,b_1,\ldots,b_r;q)_k}}
z^k,
$$
where the $q$-shifted factorials are defined by
$$
(a_1,\ldots,a_r;q)_k =\prod_{i=1}^r (a_i;q)_k,
\quad (a;q)_k=\prod_{i=0}^{k-1} (1-aq^i).
$$
The empty product equals $1$ by definition. Since $0<q<1$, $(a;q)_\infty$
is well-defined.

The other generalised matrix elements are also explicitly known in terms of
the Askey-Wilson polynomials, cf. \cite{\NoumMJASJ},
\cite{\NoumMLNM}, \cite{\KoelSIAM}, \cite{\KoelAAM}, and we need the following
special cases;
$$
\aligned
b^l_{n,0}(\infty,\s) &= C_n(\s)
\bigl( \d_{\infty,\s-1}\g_{\infty,\s} \bigr)^n
P^{(n,n)}_{l-n}(\ris;q^{2\s},1 ; q^2), \\
b^l_{0,n}(\t,\infty) &= C_n(\t)
\bigl( \d_{\t-1,\infty}\b_{\t,\infty} \bigr)^n
P^{(n,n)}_{l-n}(\rti;q^{2\t}, 1; q^2),\\
b^l_{-n,0}(\infty,\s) &= C_n(\s)
\bigl( \b_{\infty,\s-1}\a_{\infty,\s} \bigr)^n
P^{(n,n)}_{l-n}(\ris;q^{2\s+2n}, q^{2n}; q^2), \\
b^l_{0,-n}(\t,\infty) &= C_n(\t)
\bigl(\g_{\t-1,\infty}\a_{\t,\infty}  \bigr)^n
P^{(n,n)}_{l-n}(\rti;q^{2\t+2n}, q^{2n}; q^2)
\endaligned
\tag\eqname{\vglasspheleltwithoneinfty}
$$
with the constant given by
$$
C_n(\s) = q^{\hf(l-n)(l-n-1)} q^{-\s l}
{{(q^{2n+2},-q^{2\s-2l};q^2)_{l-n}}\over{
\sqrt{ (q^2,q^{2l+2n+2};q^2)_{l-n}} }}.
$$
To explain the notation in \thetag{\vglasspheleltwithoneinfty} we have the
following limit cases of $\r$;
$$
\aligned
\rti &= \lim_{\s\to\infty} 2q^{\s+\t-1}\r
= iq^{\t-1}(q\d\g+\b\a) +q^{-1}(1-q^{2\t})\b\g, \\
\ris &= \lim_{\t\to\infty} 2q^{\s+\t-1}\r
= -iq^{\s-1}(\d\b+q\g\a) +q^{-1}(1-q^{2\s})\b\g.
\endaligned
\tag\eqname{\vgllimrnaarrist}
$$
The elements $\a_{\t,\infty}$, $\a_{\infty,\s}$, etcetera, are the
straightforward limit cases of the elements defined by
$$
\aligned
\a_{\t,\s} &= q^{1/2}\a-iq^{\s-1/2}\b +iq^{\t+1/2}\g + q^{\s+\t-1/2}\d, \\
\b_{\t,\s} &= -q^{\s+1/2}\a-iq^{-1/2}\b -iq^{\s+\t+1/2}\g + q^{\t-1/2}\d, \\
\g_{\t,\s} &= -q^{\t+1/2}\a+iq^{\t+\s-1/2}\b +iq^{1/2}\g + q^{\s-1/2}\d, \\
\d_{\t,\s} &= q^{\t+\s+1/2}\a+iq^{\t-1/2}\b -iq^{\s+1/2}\g + q^{-1/2}\d.
\endaligned
\tag\eqname{\vgldefabgdst}
$$
These elements are up to a constant factor the generalised matrix elements
$b^{1/2}_{i,j}(\t,\s)$, cf. \cite{\KoelAAM, prop.~6.5}.
The $\ast$-operator on these simple generalised matrix elements is given by
$$
\a_{\t,\s}^\ast = q\d_{\t-1,\s-1},\quad
\b_{\t,\s}^\ast = -\g_{\t-1,\s+1},\quad
\g_{\t,\s}^\ast = -\b_{\t+1,\s-1},\quad
\d_{\t,\s}^\ast = q^{-1}\a_{\t+1,\s+1}.
\tag\eqname{\vglastonminelts}
$$
Finally, the polynomials involved in \thetag{\vglasspheleltwithoneinfty} are
the big $q$-Jacobi polynomials introduced by Andrews and Askey
\cite{\AndrAbqJ};
$$
P_n^{(\a,\b)}(x;c,d;q) = {}_3 \varphi_2 \left( {{q^{-n}, q^{n+\a+\b+1},
q^{1+\a}x/c}\atop{q^{\a+1},-q^{1+\a}d/c}};q,q\right).
\tag\eqname{\vglbigqJac}
$$
The special case \thetag{\vglasspheleltwithoneinfty} is obtained from the
general expression for the generalised matrix elements $b^l_{i,j}(\t,\s)$
in terms of Askey-Wilson polynomials by use of the limit
transition of the Askey-Wilson polynomials to the big $q$-Jacobi polynomials
as described by Koornwinder \cite{\KoorZSE, (6.2)}.
Equation \thetag{\vglasspheleltwithoneinfty} has been obtained by Noumi and
Mimachi \cite{\NoumMCM, thm.~3.5} in the context of a slightly different
algebra and we make a precise identification with
the algebra considered in \cite{\NoumMCM} in remark 3.5. 

\head\newsection . Basis of the representation space
\endhead

The standard orthonormal basis $\{ e_n\mid n\in\Zp\}$ of the representation
space $\Hi$ of the irreducible $\ast$-representation of $\A$ is not well
suited to calculate the action of the operators
$\p\bigl (b^l_{n,0}(\infty,\s)\bigr)$.
In this section we show that these operators act nicely in a suitable basis
of $\Hi$, which has already been introduced in \cite{\KoelAF}.

\proclaim{Proposition \theoremname{\propeigvectpris}}
{\rm (\cite{\KoelAF, prop.~4.1, cor.~4.2})}
$\Hi$ has an orthogonal basis of
eigenvectors $v_\l(q^\s)$, where $\l=-q^{2n}$, $n\in\Zp$, and
$\l = q^{2\s + 2n}$,
$n\in\Zp$, for the eigenvalue $\l$ of the self-adjoint operator $\p(\ris)$.
Moreover,
$$
\align
\langle v_\l(q^\s),v_\l(q^\s)\rangle &= q^{-2n} (q^2;q^2)_n
(-q^{2-2\s};q^2)_n (-q^{2\s};q^2)_\infty,\qquad \l=-q^{2n},\\
\langle v_\l(q^\s),v_\l(q^\s)\rangle &= q^{-2n} (q^2;q^2)_n
(-q^{2+2\s};q^2)_n (-q^{-2\s};q^2)_\infty,\qquad \l=q^{2\s+2n}.
\endalign
$$
\endproclaim

\demo{Remark \theoremname{\remsp}} Koornwinder \cite{\KoorZSE}, see also
Noumi and Mimachi \cite{\NoumMCM}, has
proved that the Haar functional on the commutative
$\ast$-subalgebra of $\A$ generated by $\ris$
is given by $h(p(\ris))=\int_{-1}^{q^{2\s}} p(x)\, d_{q^2}x$,
see \cite{\GaspR, \S 1.11} for the definition of a $q$-integral.
It is remarkable that the spectrum of the self-adjoint
operator $\p(\ris)$ equals the support of the orthogonality measure
$d_{q^2}x$ for the big $q$-Legendre polynomials on $[-1,q^{2\s}]$. This
can be used to give a different proof for this expression of the Haar
functional as a $q$-integral on $[-1,q^{2\s}]$. By determining the
spectrum of $\p(\r)$ a similar proof can be
given for the Haar functional on the commutative
$\ast$-algebra generated by $\r$ as an Askey-Wilson integral,
see \cite{\KoelV} for details. This result is due to Koornwinder
\cite{\KoorZSE}, who used the corepresentations of $\A$ to prove this
result. \enddemo

We use the convention that $v_\l(q^\s)=0$ for $\l\not=-q^{2n}$,
$\l\not=q^{2\s+2n}$. The eigenvector can be written as
$v_\l(q^\s)=\sum_{n=0}^\infty p_n(\l) e_n$, where $p_n(\l)$ is expressed
in terms of the Al-Salam--Carlitz polynomial $U^{(a)}_n$. These polynomials
can be considered as the Hermite case of the big $q$-Jacobi polynomials.
This proposition is proved as follows. The equation
$\p(\ris)\sum_{n=0}^\infty p_n(\l) e_n=\l\sum_{n=0}^\infty p_n(\l) e_n$,
$p_{-1}(\l)=0$, $p_0(\l)=1$,
leads to a three-term recurrence relation for the required polynomials
from which the identification with the Al-Salam--Carlitz polynomials can
be made. The corresponding eigenvectors are in $\Hi$ only for the given
values of $\l$, which correspond precisely to the discrete mass points
in the orthogonality measure for the Al-Salam--Carlitz polynomials.
Note that the eigenvectors $v_\l(q^\s)$ are normalised by the condition
$\langle v_\l(q^\s),e_0\rangle=p_0(\l)=1$.

The eigenvectors of $\p(\ris)$ induce an orthogonal decomposition of the
representation space $\Hi=V_1^\s\oplus V_2^\s$, where $V_1^\s$ is the
space spanned by the vectors $v_{-q^{2n}}(q^\s)$, $n\in\Zp$, and $V_2^\s$ is
the space spanned by the vectors $v_{q^{2\s+2n}}(q^\s)$, $n\in\Zp$.
The eigenvectors $v_\l(q^\s)$ of $\p(\ris)$ can be considered as a
generalisation of the standard basis, which are eigenvectors of
$\p(\rho_{\infty,\infty})=\lim_{\s\to\infty} \p(\ris)=\p(q^{-1}\b\g)$,
cf. \thetag{\vglsteropAq}. It is possible to show that the normalised
eigenvector $v_{-q^{2n}}(q^\s)/\parallel v_{-q^{2n}}(q^\s)\parallel$ tends
to $i^ne_n$ as $\s\to\infty$ and that
$v_{q^{2n+2\s}}(q^\s)/\parallel v_{q^{2n+2\s}}(q^\s)\parallel$ tends to
zero as $\s\to\infty$.

This basis of the representation space $\Hi$ is well suited for the
calculation of the action of the associated spherical elements
\thetag{\vglasspheleltwithoneinfty} in this representation. In order to show
this we go back to the Hopf $\ast$-algebra to find suitable expressions for
$\r$.

\proclaim{Proposition \theoremname{\propfacrst}}
For $\s,\t\in\R\cup\{\infty\}$ we have in $\A$
$$
\align
\b_{\t+1,\s-1}\g_{\t,\s}&=2q^{\t+\s}\r-q^{2\s-1}-q^{2\t+1}, \\
\g_{\t-1,\s+1}\b_{\t,\s}&=2q^{\t+\s}\r-q^{2\s+1}-q^{2\t-1}, \\
\a_{\t+1,\s+1}\d_{\t,\s}&=2q^{\t+\s+1}\r+1+q^{2\s+2\t+2}, \\
\d_{\t-1,\s-1}\a_{\t,\s}&=2q^{\t+\s-1}\r+1+q^{2\s+2\t-2}.
\endalign
$$
\endproclaim

\demo{Proof} This can be proved straightforwardly from
\thetag{\vgldefabgdst}, \thetag{\vgldefrhosigmatau} and the commutation
relations \thetag{\vglcommrelAq}. A more conceptual proof uses \cite{\KoelAAM,
prop.~6.5, thm.~5.1}, which imply that each product on the right hand
side is a polynomial in $\r$. The Clebsch-Gordan series, cf. e.g.
\cite{\KoelAAM, \S 4.7}, imply that this polynomial is of degree $1$.
A simple calculation determines the coefficients in this polynomial.
\qed\enddemo

\proclaim{Corollary \theoremname{\corcommrst}}
The following relations hold in $\A$;
$$
\gather
\a_{\t,\s}\r=\rho_{\t-1,\s-1}\a_{\t,\s}, \qquad
\b_{\t,\s}\r=\rho_{\t-1,\s+1}\b_{\t,\s}, \\
\g_{\t,\s}\r=\rho_{\t+1,\s-1}\g_{\t,\s}, \qquad
\d_{\t,\s}\r=\rho_{\t+1,\s+1}\d_{\t,\s}.
\endgather
$$
In particular,
$$
\gather
\a_{\infty,\s}\ris=q^2\rho_{\infty,\s-1}\a_{\infty,\s}, \qquad
\b_{\infty,\s}\ris=\rho_{\infty,\s+1}\b_{\infty,\s}, \\
\g_{\infty,\s}\ris=\rho_{\infty,\s-1}\g_{\infty,\s}, \qquad
\d_{\infty,\s}\ris=q^{-2}\rho_{\infty,\s+1}\d_{\infty,\s}.
\endgather
$$
\endproclaim

\demo{Proof} The special case follows from \thetag{\vgllimrnaarrist}. The
proof of the first statements are all similar. To prove the first, multiply
the last equation of proposition \thmref{\propfacrst} by $\a_{\t,\s}$
and use the third equation in the left hand side. Cancelling terms proves
the first statement of the corollary.
\qed\enddemo

\demo{Remark \theoremname{\remrelNM}} The elements $\a_{\infty,\s}$,
$\b_{\infty,\s}$, $\g_{\infty,\s}$ and $\d_{\infty,\s}$ are closely related
to the algebra corresponding to a quantum $3$-sphere described by Noumi and
Mimachi \cite{\NoumMCM}. To see this we calculate the commutation relations
among these elements, which can be done straightforwardly or by using
\cite{\KoelAAM, props.~6.4, 6.5} as in the proof of proposition
\thmref{\propfacrst}. We get, cf. \thetag{\vglcommrelAq},
$$
\gather
\a_{\t-1,\s+1}\b_{\t,\s} = q \b_{\t-1,\s-1}\a_{\t,\s}, \quad
\a_{\t+1,\s-1}\g_{\t,\s} = q \g_{\t-1,\s-1}\a_{\t,\s}, \\
\g_{\t+1,\s+1}\d_{\t,\s} = q \d_{\t+1,\s-1}\g_{\t,\s}, \quad
\b_{\t+1,\s+1}\d_{\t,\s} = q \d_{\t-1,\s+1}\b_{\t,\s}.
\endgather
$$
There are no easy analogues of the commutation relations for $\b$ and $\g$
and for $\a$ and $\d$, but for $\t=\infty$ we get from proposition
\thmref{\propfacrst} and the limit transition \thetag{\vgllimrnaarrist}
$$
\a_{\infty,\s+1}\d_{\infty,\s} - \d_{\infty,\s-1}\a_{\infty,\s} =
q \b_{\infty,\s-1}\g_{\infty,\s}-q^{-1}\g_{\infty,\s+1}\b_{\infty,\s}.
$$
We consider $q^\s$ as a self-adjoint element of $\A$, which commutes with $\a$,
$\b$, $\g$ and $\d$. Now we formally adjoin the unitary operator $S$ to
this algebra as in \cite{\NoumMCM, \S5}, i.e. $S$ commutes with $\a$,
$\b$, $\g$ and $\d$, and $Sq^\s=q q^\s S$, or $S$ is a shift in $\s$.
Now define $\tilde x= q^\hf S^{-1}\d_{\infty,\s}$,
$\tilde y= q^{-\hf} S\a_{\infty,\s}$, $\tilde u= q^{-\hf} S\g_{\infty,\s}$
and $\tilde v= q^\hf S^{-1}\b_{\infty,\s}$,
then these elements satisfy the relations of the quantum $3$-sphere
of Noumi and Mimachi \cite{\NoumMCM, \S 0}. The elements $c$, $d$ and $z$
defined by Noumi and Mimachi \cite{\NoumMCM, \S 0, (2.15)} correspond to
$1$, $q^{2\s}$ and $-\ris$.
\enddemo

\proclaim{Proposition \theoremname{\propactionabgd}} For $\l=-q^{2n}$,
$\l=q^{2\s+2n}$, $n\in\Zp$, we have
$$
\gather
\p(\a_{\infty,\s})v_\l(q^\s) = iq^{\hf-\s}(1+\l) v_{\l/q^2}(q^{\s-1}), \qquad
\p(\b_{\infty,\s})v_\l(q^\s) = iq^\hf (q^{2\s}-\l) v_\l(q^{\s+1}),\\
\p(\g_{\infty,\s})v_\l(q^\s) = iq^\hf v_\l(q^{\s-1}),\qquad
\p(\d_{\infty,\s})v_\l(q^\s) = -iq^{\hf+\s} v_{\l q^2}(q^{\s+1}).
\endgather
$$
\endproclaim

\demo{Proof} First observe that
the result for $\p(\g_{\infty,\s})$, respectively $\p(\d_{\infty,\s})$,
implies the result for $\p(\b_{\infty,\s})$, respectively
$\p(\g_{\infty,\s})$, by propositions \thmref{\propfacrst} and
\thmref{\propeigvectpris}.

Consider $\p(\g_{\infty,\s})$, then we see that
$\p(\g_{\infty,\s})v_\l(q^\s)=Cv_\l(q^{\s-1})$ for some constant $C$
by corollary \thmref{\corcommrst} and proposition \thmref{\propeigvectpris}.
To calculate $C$ we use $\langle v_\l(q^\s),e_0\rangle = 1$ and so
$$
C=\langle \p(\g_{\infty,\s})v_\l(q^\s),e_0\rangle
=\langle v_\l(q^\s),-\p(\b_{\infty,\s-1})e_0\rangle = iq^\hf
$$
by \thetag{\vglastonminelts}, \thetag{\vgldefabgdst} and
\thetag{\vgldefreprA}.

Similarly, $\p(\d_{\infty,\s})v_\l(q^\s)=Cv_{\l q^2}(q^{\s+1})$ for some
constant $C$ with
$$
C=\langle \p(\d_{\infty,\s})v_\l(q^\s),e_0\rangle
=\langle v_\l(q^\s),q^{-1}\p(\a_{\infty,\s+1})e_0\rangle = -iq^{\s+\hf},
$$
which proves the proposition.
\qed\enddemo

\demo{Remark \theoremname{\remcqrel}} If we use
$v_\l(q^\s)=\sum_{n=0}^\infty p_n(\l)e_n$ and we work
out what proposition \thmref{\propactionabgd} means for the Al-Salam--Carlitz
polynomials, which can be expressed as a ${}_2\vp_1$-series, then we see that
the expressions in proposition \thmref{\propactionabgd} are equivalent to
one of Heine's contiguous relations for the ${}_2\vp_1$-series, cf.
\cite{\GaspR, Ex.~1.9(ii)}. \enddemo

Since $-q\p(\g)=\p(\b)$ implies $\p(\r)=\p(\rho_{\s,\t})$,
$\p(\a_{\t,\s})=\p(\a_{\s,\t})$,
$\p(\b_{\t,\s})=\p(\g_{\s,\t})$,
$\p(\g_{\t,\s})=\p(\b_{\s,\t})$ and
$\p(\d_{\t,\s})=\p(\d_{\s,\t})$,
we obtain the following proposition as a corollary of propositions
\thmref{\propeigvectpris} and \thmref{\propactionabgd}.

\proclaim{Proposition \theoremname{\coreigvectprti}}
$\Hi$ has an orthogonal basis of
eigenvectors $v_\l(q^\t)$, where $\l=-q^{2n}$, $n\in\Zp$, $\l=q^{2\t+2n}$,
$n\in\Zp$, for the eigenvalue $\l$ of the self-adjoint operator $\p(\rti)$.
For $\l=-q^{2n}$, $\l=q^{2\t+2n}$, $n\in\Zp$, we have
$$
\gather
\p(\a_{\t,\infty})v_\l(q^\t) = iq^{\hf-\t}(1+\l) v_{\l/q^2}(q^{\t-1}), \qquad
\p(\b_{\t,\infty})v_\l(q^\t) = iq^\hf v_\l (q^{\t-1}),\\
\p(\g_{\t,\infty})v_\l(q^\t) = iq^\hf (q^{2\t}-\l)v_{\l}(q^{\t+1}),\qquad
\p(\d_{\t,\infty})v_\l(q^\t) = -iq^{\hf+\t} v_{\l q^2}(q^{\t+1}).
\endgather
$$
\endproclaim

It follows from propositions \thmref{\propactionabgd},
\thmref{\coreigvectprti} and \thetag{\vglasspheleltwithoneinfty}
that for $n\geq 0$
$$
\aligned
\p(b^l_{n,0}(\infty,\s)) v_\l(q^\s) &=
C_n(\s) q^{n\s} P^{(n,n)}_{l-n}(\l;q^{2\s},1 ; q^2)\, v_{\l q^{2n}}(q^\s),\\
\p(D.b^l_{0,n}(\t,\infty))  v_\m(q^\t) &=
C_n(\t) q^{n(\t+1)} P^{(n,n)}_{l-n}(\m;q^{2\t},1 ; q^2)
\, v_{\m q^{2n}}(q^\t)
\endaligned
\tag\eqname{\vglsombJpos}
$$
and
$$
\aligned
\p(b^l_{-n,0}(\infty,\s)) v_\l(q^\s) &=
C_n(\s)(-1)^n q^{n(\s-1)} (-\l,\l q^{-2\s};q^{-2})_n \\
&\qquad\qquad\qquad\times
P^{(n,n)}_{l-n}(\l q^{-2n};q^{2\s},1 ; q^2)\, v_{\l q^{-2n}}(q^\s),\\
\p(D.b^l_{0,-n}(\t,\infty)) v_\m(q^\t) &=
C_n(\t)(-1)^n q^{n(\t-2)} (-\m,\m q^{-2\t};q^{-2})_n\\
&\qquad\qquad\qquad\times
P^{(n,n)}_{l-n}(\m q^{-2n};q^{2\t},1 ; q^2)\, v_{\m q^{-2n}}(q^\t),
\endaligned
\tag\eqname{\vglsombJneg}
$$
where we used $D.\rti=\rti$, $D.\a_{\t,\infty}=q^{-\hf}\a_{\t,\infty}$,
$D.\b_{\t,\infty}=q^\hf\b_{\t,\infty}$,
$D.\g_{\t,\infty}=q^{-\hf}\g_{\t,\infty}$ and
$D.\d_{\t,\infty}=q^\hf\d_{\t,\infty}$ and
$P^{(\a,\b)}_n(Ax;Ac,AD;q)=P^{(\a,\b)}_n(x;c,D;q)$ for $A>0$.

From \thetag{\vglsombJpos} and \thetag{\vglsombJneg} it follows that the
representation operators corresponding to the right hand side of
\thetag{\vglabstractadfeen} preserve the orthogonal decomposition
$\Hi\hat\otimes\Hi=\bigoplus_{i,j=1}^2 V_i^\t\hat\otimes V_j^\s$.

\demo{Remark \theoremname{\remgensit}} It is also possible to determine
the spectrum of the
self-adjoint operator $\p(\r)$ and to use generalised eigenvectors of
$\p(\r)$ in the representation space $\Hi$. It is then possible to prove
an analogue of proposition \thmref{\propactionabgd} and to calculate the
action of the associated spherical elements as in
\thetag{\vglsombJpos} and \thetag{\vglsombJneg}, but it turns out that
only for the choices of the parameters as in this section
the orthogonal basis is preserved under the action of the associated
spherical elements. This means that restriction $\mu=\infty$ in
\thetag{\vglabstractadfeen} is necessary in order to use the method of
this paper to convert \thetag{\vglabstractadfeen} into an
addition formula in commuting variables.
\enddemo

\head\newsection . Addition formula
\endhead

In this section we derive an explicit addition formula for a two-parameter
set of Askey-Wilson polynomials. We show that two known addition formulas
can be obtained as limit cases.

From \thetag{\vglsombJpos} and \thetag{\vglsombJneg} we see how the right
hand side of \thetag{\vglabstractadfeen} acts in the
basis $v_\m(q^\t)\otimes v_\l(q^\s)$ of the representation space
$\Hi\hat\otimes\Hi$. So we now consider the left hand side of
\thetag{\vglabstractadfeen}.
Since $\Delta$ is a $\ast$-algebra homomorphism, the
left hand side of \thetag{\vglabstractadfeen} is a polynomial in
$\Delta(\r)$. Now we study the self-adjoint operator
$(\p\otimes\p)\Delta(\r)$ in $\Hi\hat\otimes\Hi$. First we consider
$$
\multline
\Delta(\b_{\t+1,\s-1}\g_{\t,\s}) =
(D.\a_{\t+1,\infty})(D.\g_{\t,\infty})\otimes\b_{\infty,\s-1}\a_{\infty,\s}
+\\
(D.\a_{\t+1,\infty})(D.\d_{\t,\infty})\otimes
\b_{\infty,\s-1}\g_{\infty,\s}
+(D.\b_{\t+1,\infty})(D.\g_{\t,\infty})\otimes
\d_{\infty,\s-1}\a_{\infty,\s}\\
+(D.\b_{\t+1,\infty})(D.\d_{\t,\infty})\otimes\d_{\infty,\s-1}\g_{\infty,\s},
\endmultline
\tag\eqname{\vglDbg}
$$
because we can apply \thetag{\vglabstractaddf} for $l=\hf$, $\m=\infty$ twice,
since the elements in \thetag{\vgldefabgdst} are generalised matrix elements.
We obtain, cf. propositions \thmref{\propfacrst}, \thmref{\propactionabgd}
and \thmref{\coreigvectprti}
$$
\gather
2(\p\otimes\p)\Delta(\r)v_\m(q^\t)\otimes v_\l(q^\s) =
\tag\eqname{\vgldtrecDrts} \\
 q^2v_{\m q^2}(q^\t)\otimes v_{\l q^2}(q^\s)
+ q^{-2}(1+\l)(1+\m)(1-\l q^{-2\s})(1-\m q^{-2\t})
v_{\m q^{-2}}(q^\t)\otimes v_{\l q^{-2}}(q^\s)\\
+ \Bigl( \l q^{1-\s}(q^{-\t}-q^\t)+\m q^{1-\t}(q^{-\s}-q^\s) +\l\m
q^{1-\t-\s}(1+q^2)\Bigr) v_\m(q^\t)\otimes v_\l(q^\s).
\endgather
$$
It follows from \thetag{\vgldtrecDrts} that $(\p\otimes\p)\Delta(\r)$
preserves the orthogonal decomposition of the representation space
$\Hi\hat\otimes\Hi=\bigoplus_{i,j=1}^2 V_i^\t\hat\otimes V_j^\s$.
Moreover, it leaves the space
$\sum_{m=0}^\infty \C v_{\m q^{2m}}(q^\t)\otimes v_{\l
q^{2m}}(q^\s)$ invariant whenever $\m\in\{ -1,q^{2\t}\}$ or
$\l\in\{ -1,q^{2\s}\}$. In each $V_i^\t\otimes V_j^\s$ we have two of
such spaces, say $W_{i,j}^\t$ and $W_{i,j}^\s$, where
$W_{i,j}^\t$ is the subspace corresponding to $\m\in\{ -1,q^{2\t}\}$.
In each $W_{i,j}^\n$, $i,j=1,2$, $\n=\s,\t$, there are
$(\p\otimes\p)\Delta(\r)$-invariant subspaces. We consider the subspace
$W_{1,1}^{\t,p} =
\sum_{m=0}^\infty \C v_{-q^{2m}}(q^\t)\otimes v_{-q^{2m+2p}}(q^\s)$,
which is invariant under $(\p\otimes\p)\Delta(\r)$ for every $p\in\Zp$.

So we define $w_m$ to be the normalised vector
$v_{-q^{2m}}(q^\t)\otimes v_{-q^{2m+2p}}(q^\s)$, i.e. by propositions
\thmref{\propeigvectpris} and \thmref{\coreigvectprti},
$$
w_m = {{q^{2m+p}\ v_{-q^{2m}}(q^\t)\otimes v_{-q^{2m+2p}}(q^\s)}\over
{\sqrt{(q^2,-q^{2-2\t};q^2)_m(q^2,-q^{2-2\s};q^2)_{m+p}
(-q^{2\t},-q^{2\s};q^2)_\infty}}}.
\tag\eqname{\vglwm}
$$
Then we find the following expression for $(\p\otimes\p)\Delta(\r)$;
$$
\gathered
2(\p\otimes\p)\Delta(\r) w_m = a_m w_{m+1} + b_m w_m + a_{m-1} w_{m-1}, \\
a_m = \sqrt{ (1-q^{2m+2})(1-q^{2m+2p+2})(1+q^{2-2\t+2m})
(1+q^{2+2p-2\s+2m}) }, \\
b_m = q^{2m} \Bigr( q^{1+2p-\s}(q^\t-q^{-\t})+ q^{1-\t}(q^\s-q^{-\s}) +
q^{2m+2p+1-\t-\s}(1+q^2)\Bigr).
\endgathered
\tag\eqname{\vglthreetermwm}
$$

The three-term recurrence for orthonormal polynomials corresponding to
\thetag{\vglthreetermwm} can be solved in terms of Askey-Wilson polynomials,
or in terms of the $q$-Laguerre case of the $q$-Jacobi polynomials defined in
\thetag{\vgldefAWalsqJacobi}. So we define
$$
l_n^{(\a)}(x;s,t\mid q) = p_n(x;q^{1/2+\a}s/t,q^{1/2}t/s,
-q^{1/2}/(st), 0\mid q),
\tag\eqname{\vgldefAWalsqLaguerre}
$$
i.e. by letting $\b\to\infty$ in \thetag{\vgldefAWalsqJacobi}. From the
three-term recurrence relation for the Askey-Wilson polynomials, cf. Askey and
Wilson \cite{\AskeW, \S 1}, we obtain
$$
\gather
2xl_n(x)=l_{n+1}(x) +
(1-q^n)(1-q^{\a+n})(1+q^ns^{-2})(1+q^{n+\a}t^{-2}) l_{n-1}(x)\\
+ q^n\Bigl( (t-t^{-1})q^\hf s^{-1} +
(s-s^{-1})q^{\hf+\a}t^{-1} +(1+q)q^{\hf+n+\a}s^{-1}t^{-1}\Bigr)l_n(x),
\endgather
$$
where $l_n(x)=l_n^{(\a)}(x;s,t\mid q)$. The orthogonality measure for the
$q$-Laguerre polynomials follows from \cite{\AskeW, thm.~2.4}, and as
in \cite{\KoorZSE, (2.7)} we denote the normalised orthogonality measure
by $dm^{(\a)}(\cdot;s,t\mid q)$.

If we denote the corresponding
orthonormal polynomials by $\tilde l_n^{(\a)}(x;s,t\mid q)$, we see
that the three-term recurrence \thetag{\vglthreetermwm}
is solved by the $q$-Laguerre
polynomials $\tilde l_m^{(p)}(\cdot;q^\t,q^\s\mid q^2)$.
The spectral theory of Jacobi matrices shows that the spectral
decomposition $E$ of the self-adjoint operator $(\p\otimes\p)\Delta(\r)$
on the Hilbert space $W_{1,1}^{\t,p}$  with orthonormal basis $w_m$,
$m\in\Zp$, is given by
$$
\langle E(B) w_n, w_m\rangle = \int_B \bigl( \tilde l_n^{(p)}
\tilde l_m^{(p)}\bigr) (x;q^\t,q^\s\mid q^2) \,
dm^{(p)}(x;q^\t,q^\s\mid q^2)
$$
for $B$ a Borel subset of $\R$, cf. Berezanski\u\i\
\cite{\Bere, Ch.~VII,  \S 1}, Dombrowski \cite{\Domb}.
The mapping
$\Lambda\colon W_{1,1}^{\t,p}\to L^2(\R,dm^{(p)}(x;q^\t,q^\s\mid q^2))$
given by
$\Lambda(w_m)(x)=\tilde l_m^{(p)}(x;q^\t,q^\s\mid q^2)$ extends to a unitary
mapping, since the corresponding moment problem is determined. $\Lambda$
intertwines $(\p\otimes\p)\Delta(\r)$ on $W_{1,1}^{\t,p}$
with the multiplication operator $M$, $Mf(x)=xf(x)$, on this $L^2$-space.

Next we apply $\Lambda\circ p_l^{(0,0)}\bigl(
(\p\otimes\p)\Delta(\r);q^\t,q^\s\mid q^2\bigr)$
to $w_m$ given by \thetag{\vglwm} on the one hand. On the other hand we apply
$\p\otimes\p$ of the right hand side of \thetag{\vglabstractadfeen}
for $\m=\infty$ to the vector $w_m$ for which we use
\thetag{\vglsombJpos} and \thetag{\vglsombJneg} before we apply
$\Lambda$. This results in the addition formula for the
two-parameter $q$-Legendre polynomial $p^{(0,0)}_l(x;q^\s,q^\t\mid
q^2)$ as an identity in $L^2(\R,dm^{(p)}(x;q^\t,q^\s\mid q^2))$, but
since it only involves polynomials it holds everywhere. The following
theorem is then obtained after replacing $q^2$ by $q$ and $\t$, $\s$ by
$2\t$, $2\s$.

\proclaim{Theorem \theoremname{\thmfirstaddform}}
The following addition
formula holds for $l,m,p\in\Zp$, $\s,\t\in\R$, $x\in\C$;
$$
\align
&p^{(0,0)}_l(x;q^\t,q^\s\mid q) l_m^{(p)}(x;q^\t,q^\s\mid q) = \\
& \sum_{n=0}^l D^{n,l}(\t,\s)
P^{(n,n)}_{l-n}(-q^m;q^{2\t},1;q)
P^{(n,n)}_{l-n}(-q^{m+p};q^{2\s},1;q)
l_{m+n}^{(p)}(x;q^\t,q^\s\mid q)  \\
&+ \sum_{n=1}^l D^{n,l}(\t,\s)
(q^m,q^{m+p},-q^{m+p-2\s},-q^{m-2\t};q^{-1})_n \\
&\qquad\times P^{(n,n)}_{l-n}(-q^{m-n};q^{2\t},1;q)
P^{(n,n)}_{l-n}(-q^{m+p-n};q^{2\s},1;q^2)
l_{m-n}^{(p)}(x;q^\t,q^\s\mid q),
\endalign
$$
with the constant given by
$$
D^{n,l}(\t,\s) = (-q^{2\s-l},-q^{2\t-l};q)_{l-n}
{{(q^{l-n+1};q)_n}\over{(q;q)_n}} (q^{n+1};q)_l
q^{\hf (l-n)(l-n-2\s-2\t)}
$$
and the notation for the $q$-Legendre polynomials
$p_l^{(0,0)}(\cdot ;q^\t,q^\s\mid q)$,
the $q$-Laguerre polynomials $l_m^{(p)}(\cdot;q^\t,q^\s\mid q)$ and the big
$q$-ultraspherical polynomials $P^{(n,n)}_{l-n}(\cdot;q^{2\t},1;q)$ defined
by \thetag{\vgldefAWalsqJacobi}, \thetag{\vgldefAWalsqLaguerre} and
\thetag{\vglbigqJac}.
\endproclaim

\demo{Proof} It remains to calculate the constants involved.
If $h_m$ denotes the square norm of $l_m^{(p)}(\cdot;q^\t,q^\s\mid q^2)$,
then we use, cf. Askey and Wilson \cite{\AskeW, thm.~2.3},
$$
{{\parallel v_{-q^{2m+2n}}(q^\t)\otimes v_{-q^{2m+2n+2p}}(q^\s)\parallel}
\over{\parallel v_{-q^{2m}}(q^\t)\otimes v_{-q^{2m+2p}}(q^\s)\parallel}}
\sqrt{ {{h_m}\over{h_{m+n}}} } = q^{-2n}
$$
to see that before changing $q^2$ to $q$ and $\t$, $\s$ to
$2\t$, $2\s$ we have
$$
D^{n,l}(\t,\s) = q^l (q^{2l+2};q^2)_l C_n(\s)C_n(\t)
q^{n(\s+\t-1)}.
$$
A straightforward calculation yields the result.
\qed\enddemo

\demo{Remark \theoremname{\remsymmaddform}}
In deriving theorem \thmref{\thmfirstaddform} we have
chosen a certain
generalised eigenvector of $(\p\otimes\p)\Delta(\r)$ from the
space $W_{1,1}^{\t,p}\subset\Hi\hat\otimes\Hi$. It is also possible to chose
a generalised eigenvector from similarly defined spaces
$W_{i,j}^{\n,p}$, $i,j=1,2$, $\n=\s,\t$, so that we obtain
in total eight of such eigenvectors. This results in eight of these addition
formulas, which are obviously four by four equivalent by interchanging $\t$
and $\s$ and noting that the $q$-Legendre polynomial involved is also
invariant under such a change. The remaining four types of addition formulas
differ, since the big $q$-ultraspherical polynomials are evaluated at other
points of the spectrum. However, from
$$
P^{(\a,\b)}_n (-x;c,d;q) = (-q^{\a-\b}d/c)^n
{{(q^{\b+1},-q^{\b+1}c/d;q)_n}\over{(q^{\a+1},-q^{\a+1}d/c;q)_n}}
P^{(\b,\a)}_n (x;d,c;q),
\tag\eqname{\vglbigqJacflip}
$$
which is a direct consequence of the orthogonality relations, cf.
\cite{\AndrAbqJ}, and the trivial relation
$P^{(\a,\b)}_n (Ax;Ac,Ad;q) = P^{(\a,\b)}_n (x;c,d;q)$, for $A>0$,
we see that
$$
P_{l-n}^{(n,n)}(-q^m;q^{2\t},1;q) = (-q^{-2\t})^{l-n}
{{(-q^{n+1+2\t};q)_{l-n}}\over{(-q^{n+1-2\t};q)_{l-n}}}
P_{l-n}^{(n,n)}(q^{m-2\t};q^{-2\t},1;q).
\tag\eqname{\vglswitchpoints}
$$
If we use this in theorem \thmref{\thmfirstaddform} and we next change $x$
into $-x$ using
$p_n^{(\a,\b)}(-x;s,t\mid q)=(-1)^n p_n^{(\a,\b)}(x;-s,t\mid q)$
for the $q$-Legendre polynomial, $\a=\b=0$, and for the $q$-Laguerre
polynomial, $\b\to\infty$, and we change
$q^\t$ to $-q^{-\t}$, then we obtain the
same addition formula as if we had started off with the space
$W_{2,1}^{\t,p}$, i.e. with vectors of the type
$v_{q^{2\t+2m}}(q^\t)\otimes v_{-q^{2m+2p}}(q^\s)$, $m\in\Zp$.
A similar approach can be used for the other big $q$-ultraspherical
polynomial, and a combination of both shows that theorem
\thmref{\thmfirstaddform} contains all the other possibilities by symmetry
considerations.

We also obtain the same addition formula if we would have started with one of
the other irreducible $\ast$-representations of the Hopf $\ast$-algebra $\A$.
\enddemo

\demo{Remark \theoremname{\remspecialone}} Theorem \thmref{\thmfirstaddform}
contains as special
cases two addition formulas for $q$-Legendre polynomials previously obtained
from the quantum $SU(2)$ group. Firstly, we rescale $x=q^{1-\s-\t}y/2$
and we take the limits $\t\to\infty$. The limit transition of
Askey-Wilson polynomials to big $q$-Jacobi polynomials as in
\cite{\KoorZSE, (6.2)} shows that
$p^{(0,0)}_l(x;q^\t,q^\s\mid q)$ and $l_m^{(p)}(x;q^\t,q^\s\mid q)$
tend to the big $q$-Legendre and big $q$-Laguerre polynomials. Moreover,
the big $q$-ultraspherical polynomial $P^{(n,n)}_{l-n}(x;q^{2\t},1;q)$
tends to a little $q$-Jacobi polynomial, and we obtain an addition
formula for a big $q$-Legendre polynomial. If we next let the order $p$
of the big $q$-Laguerre polynomial tend to infinity, so that the big
$q$-Laguerre polynomials tend to the Al-Salam--Carlitz $U_n^{(a)}$, we
obtain the addition formula for the big $q$-Legendre polynomials derived
in \cite{\KoelAF, thm.~4.3}. For this we have to observe that a big
$q$-ultraspherical polynomial of argument zero can be rewritten as a dual
$q$-Krawtchouk polynomial by \cite{\GaspR, (III.5), (III.6)}.

Secondly, rescaling $x=q^{1-\s-\t}y/2$ and letting both $\s,\t\to\infty$,
and using the limit transition of the Askey-Wilson polynomials to the
little $q$-Jacobi polynomial, cf. \cite{\KoorZSE, (6.4)}, and the limit
transition of the big $q$-Jacobi polynomials to the little $q$-Jacobi
polynomials, shows that the addition formula in theorem
\thmref{\thmfirstaddform} contains Koornwinder's addition formula for
the little $q$-Legendre polynomials \cite{\KoorAF, thm.~4.1}
as a special case, see Rahman \cite{\Rahm} for an analytic proof.
\enddemo

\demo{Remark \theoremname{\remcontqLeg}}
Another special case of theorem \thmref{\thmfirstaddform} is the case
$\s=\t=0$, which gives an addition formula for the continuous $q$-Legendre
polynomial. This addition formula is an alternative for the addition formula
for the continuous $q$-Legendre polynomials derived by Rahman and Verma
\cite{\RahmV, (1.24) with $a=q^{1/4}$}, see also \cite{\KoelSIAM, \S 4}.
\enddemo

\head\newsection . Product formula
\endhead

Using the orthogonality relations for the $q$-Laguerre polynomial a
product formula can be easily obtained from the addition formula.
For this we need the normalised orthogonality measure
$dm^{(\a)}(x;s,t\mid q)$ for the $q$-Laguerre polynomials;
$$
\int_\R \Bigl( l^{\a)}_m l^{\a)}_n\Bigr) (x;s,t\mid q)
\, dm^{(\a)}(x;s,t\mid q) =\d_{n,m}
(q,q^{1+\a},-q^{-2},-q^{1+\a}t^{-2}:q)_n,
\tag\eqname{\vglorthoqLag}
$$
cf. Askey and Wilson \cite{\AskeW, thm.~2.5} in combination with
\thetag{\vgldefAWalsqLaguerre} for the explicit form of the orthogonality
measure as an absolutely continuous measure on $[-1,1]$ plus a
finite number, possibly zero, of discrete mass points off $[-1,1]$.

Now multiply theorem \thmref{\thmfirstaddform} by
$l^{(p)}_{m+r}(x;q^\t,q^\s\mid q)$ and integrate against the normalised
orthogonality measure $dm^{(p)}(x;q^\t,q^\s\mid q^2)$
use the orthogonality relations \thetag{\vglorthoqLag}
for the $q$-Laguerre polynomials to obtain
the following product formula for the product of two big $q$-ultraspherical
polynomials, which for $n=0$ gives a $q$-analogue of the product formula
for the Legendre polynomials \thetag{\vglprodformultra}.

\proclaim{Corollary \theoremname{\corprodformeen}} For $l,n\in\Zp$,
$0\leq n\leq l$ we have the product formula
$$
\multline
P^{(n,n)}_{l-n}(-q^m;q^{2\t},1;q)
P^{(n,n)}_{l-n}(-q^{m+p};q^{2\s},1;q) = \\{1\over C}
\int_\R p^{(0,0)}_l(x;q^\t,q^\s\mid q)
l_m^{(p)}(x;q^\t,q^\s\mid q)
l_{m+n}^{(p)}(x;q^\t,q^\s\mid q)\, dm^{(p)}(x;q^\t,q^\s\mid q) ,
\endmultline
$$
with
$$
C = D^{n,l}(\t,\s) (q,q^{1+p},-q^{1-2\t},-q^{1+p-2\s};q)_{m+n}.
$$
\endproclaim

The same product formula is obtained if we multiply theorem
by $l^{(p)}_{m-r}(x;q^\t,q^\t\mid q)$ before
using the orthogonality relations.
A product formula for big $q$-ultraspherical polynomials at other points of
the spectrum follows from \thetag{\vglswitchpoints}.

\head\newsection . The limit case $q\uparrow 1$
\endhead

We consider the limit transition $q\uparrow 1$ of the
addition and product formulas from theorem \thmref{\thmfirstaddform} and
corollary \thmref{\corprodformeen} to the addition and product formulas for
the Legendre polynomials as presented in \S 1.
This can be done using the theorems of Van Assche
and Koornwinder \cite{\VAsscK, thm.~1,~2}, cf. \cite{\KoelAF, \S 5} for
another example of the use of these theorems.

We divide both sides of the addition formula in theorem
\thmref{\thmfirstaddform} by $l^{(p)}_m(x;q^\t,q^\s;\mid q) (q;q)_l$
and we substitute $q=c^{1/m}$, $0<c<1$, so that $q\uparrow 1$ corresponds
to $m\to\infty$, and $p=mr$ with $r= \ln p/\ln c\in\Zp$. Also, put $t=c^r$.
Under these identifications we see that the various $q$-Jacobi polynomials
tend to the Jacobi polynomials as $m\to\infty$, e.g.
$$
\gather
P_{l-n}^{(n,n)}(-q^{m+p-n};q^{2\s},1;q)\to R_{l-n}^{(n,n)}(-ct),\\
(q;q)_l^{-1} p_l^{(0,0)}(x;q^\t,q^\s\mid q) \to
4^l R_l^{(0,0)}(x).
\endgather
$$
Also
$$
\gather
{{D^{n,l}(\t,\s)}\over{(q;q)_l}} \to 4^{l-n} {{(l-n+1)_n
(n+1)_l}\over{n!\, l!}}, \\
(q^m,q^{m+p},-q^{m+p-2\s},-q^{m-2\t};q^{-1})_n \to
\bigl( (1-c^2)(1-c^2t^2)\bigr)^n.
\endgather
$$
The non-trivial limit transition of the quotient of two $q$-Laguerre
polynomials is handled using theorem~1 of Van Assche and Koornwinder
\cite{\VAsscK} and we obtain for $k\in\Z$
$$
\lim_{m\to\infty} {{ l^{(mr)}_{m+k}(x;c^{\t/m},c^{\s/m} \mid c^{1/m})}\over
{l^{(mr)}_m(x;c^{\t/m},c^{\s/m}\mid c^{1/m})}} =
\bigl( (1-c^2)(1-t^2c^2)\bigr)^{k/2} \rho^k\bigl(
(x-c^2t)/\sqrt{(1-c^2)(1-t^2c^2)}\bigr),
$$
where $\rho(x)=x+\sqrt{x^2-1}$ ($\vert\rho(x)\vert>1$ for $x\in\C\backslash
[-1,1]$). Using these limit transitions in theorem \thmref{\thmfirstaddform}
and $2T_n(x)=\rho^k(x)+\rho^{-k}(x)$
shows that for $q\uparrow 1$ we obtain the addition formula
\thetag{\vgladdformLeg} for the Legendre polynomials.

The limit transition of the product formula of corollary
\thmref{\corprodformeen} to the product formula
\thetag{\vglprodformultra} is treated similarly, but now we have to use
theorem~2 of Van Assche and Koornwinder \cite{\VAsscK} on the weak
asymptotics of the $q$-Laguerre polynomials. A straightforward calculation
shows that this theorem implies for continuous $f$ and $n\in\Zp$
$$
\multline
\lim_{n\to\infty} \int_{\R} f(x) \biggl( \tilde l^{(mr)}_m
\tilde l^{(mr)}_{m+n}\biggr) (x;q^\t,q^\s\mid q) \,
dm^{(mr)}(x;q^\t,q^\s\mid q) = \\ {1\over \pi}
\int_{B-2A}^{B+2A} f(x)
{{T_n\bigl((x-B)/(2A)\bigr)}\over{\sqrt{4A^2-(x-B)^2}}}\, dx
\endmultline
$$
with $B=c^2t$, $2A=\sqrt{(1-c^2)(1-c^2t^2)}$. Using the previous limit
transitions of this section we see that the product formula of corollary
\thmref{\corprodformeen} tends to the classical product formula
\thetag{\vglprodformultra} as $q\uparrow 1$, see \cite{\VAsscK, \S 4},
\cite{\KoelAF, \S 5} for similar limit transitions.


\enddocument